\documentclass[11pt]{amsart}
\usepackage{amssymb, amscd, epsfig, times}

\newcommand\Lie{\operatorname{Lie}}
\newcommand\GL{\operatorname{GL}}

\newcommand\Hom{\operatorname{Hom}}
\newcommand\sym{\operatorname{Sym}}

\newcommand\Pev{\operatorname{Pev}}

\def\half{{\tfrac{1}{2}}} 
\def\kbold{{\mathbf k}} 
 
\newcommand{\Z}{\mathbb Z}
\newcommand{\C}{\mathbb C}
\newcommand{\R}{\mathbb R}
\newcommand{\D}{\mathbb D}
\newcommand{\N}{\mathbb N}
\newcommand{\mP}{\mathbb P}
\newcommand{\mB}{\mathbb B}
\numberwithin{equation}{section}

\newtheorem{theorem}{Theorem}[section]
\newtheorem{lemma}[theorem]{Lemma}
\newtheorem{definition}[theorem]{Definition}
\newtheorem{proposition}[theorem]{Proposition}

\theoremstyle{remark}
\newtheorem{remark}[theorem]{Remark}

\newtheorem{question}[theorem]{Question}

\newtheorem{example}[theorem]{Example}
\begin{document}
\title{On the geometry of the Calogero-Moser system}

\author{Wim Couwenberg}
\address{Korteweg-de Vries Instituut voor Wiskunde, 
Plantage Muidergracht 24, NL-1018 TV Amsterdam,
Nederland}\email{wcouwenb@science.uva.nl}
\author{Gert Heckman}
\address{Mathematical Institute, University of
Nijmegen, P.O. Box 9010, NL-6500 GL Nijmegen, Nederland}
\email{heckman@math.ru.nl}
\author{Eduard Looijenga}
\address{Faculteit Wiskunde en Informatica, Universiteit 
Utrecht, P.O. Box 80.010, NL-3508 TA Utrecht, Nederland} 
\email{looijeng@math.uu.nl}

\dedicatory{To Gerrit van Dijk for his 65$^{th}$ birthday}
\begin{abstract}
We discuss a special eigenstate of the quantized periodic Calogero-Moser
system associated to a root system. This state has 
the property that its eigenfunctions, when regarded as multivalued functions 
on the space of regular conjugacy classes in the corresponding semisimple
complex Lie group, transform under monodromy according to the 
complex reflection representation of the affine Hecke algebra. 
We show that this endows the space of conjugacy classes in question with a 
projective structure. For a certain parameter range this projective structure
underlies a complex hyperbolic structure. If in addition 
a Schwarz type of integrality condition is satisfied, then it 
even has the structure of a ball quotient minus a Heegner divisor. 
For example, the case of
the root system $E_8$ with the triflection monodromy representation describes
a special eigenstate for the system of 12 unordered points on the projective
line under a particular constraint.
\end{abstract}

\maketitle

\section{Introduction}
The classical Calogero-Moser system describes a finite number of
identical point particles on the real line under the influence 
of an inverse square potential. 
The relative positions of these points are parametrized by the complement
of a hyperplane arrangement, namely by $\R ^{n+1}$ modulo its main 
diagonal minus the union of the hyperplanes for which two coordinates are equal
(in case they are $n+1$ in number). 
Thus we find ourselves immediately dealing with a root system of type $A_n$.
It is therefore not surprising that there is a generalization to arbitrary finite root systems.
This can be (and has been) been taken a step further by letting the basic input be
an abstract root system $R$ and a one dimensional connected
complex Lie group $\mathbb G$: if $P$ denotes the weight lattice of $R$, then
the substitute for the configuration space is the abelian complex Lie group
$\Hom (P, \mathbb G)$ deprived from the fixed point divisors of the reflections in 
the Weyl group $W$ of $R$. 
So if $\mathbb G$ is the additive group $\C$, then we recover the previous 
generalization, but if $\mathbb G$ is the multiplicative group $\C^\times$ resp.\ 
an elliptic curve, then we get the complement of a reflection arrangement 
in an algebraic torus resp.\ an abelian variety. It is common terminology 
to talk about the non-periodic, the periodic and the doubly periodic case
of the root system generalization of the Calogero-Moser system.

\medskip\noindent
The general integrability properties of these systems are nowadays rather well 
understood. 
Moser showed that the classical ($A_n$) system is completely integrable in the
additive and multiplicative case by realizing such a system as a 
Lax pair \cite{M75}. Olshanetsky and Perelomov extended his method to 
the elliptic curve case \cite{OP76} and also began the study of the Calogero-Moser system
in the general context of root systems \cite{OP81},\cite{OP83}.
The question of complete integrability for the quantization of the
Calogero-Moser system in the multiplicative case and for a general root system $R$
was solved by complex analytic methods in the work of Heckman and Opdam \cite{HO87/88}. 
The associated eigenvalue system is also called the hypergeometric
system associated with $R$ because in the rank one situation it 
boils down to the Euler-Gauss hypergeometric equation.
Heckman later showed that these results could be obtained in an 
elementary way by means of the theory of Dunkl operators \cite{H91}.
Subsequently Cherednik and Matsuo noticed the analogy between this approach
and the Knizhnik-Zamolodchikov equation of conformal field theory \cite{C91}, 
\cite{M92}. 
Cherednik also made the important observation that the relations for the 
(trigonometric) Dunkl operators are those of the degenerate affine Hecke 
algebra of Drinfeld and Lusztig \cite{D86}, \cite{O95}, \cite{H97}. The affine
Hecke algebra appeared in the monodromy of the hypergeometric system
associated with $R$ as well \cite{HO87/88}, and the role of the (double) affine Hecke algebra
in the basic hypergeometric extension (where $q$-difference operators
replace differential operators) became increasingly important \cite{M03}.

\medskip\noindent
It turns out that our previous work on the geometric 
structures on complements of projectivized hyperplane arrangements  
can be understood in and generalized to this context. Our paper \cite{CHL03}
accomplished that in the additive setting by generalizing the work of 
Deligne-Mostow on the Lauricella functions (which corresponds from our 
perspective to the classical $A_n$-case) to arbitrary root systems.

In the present article we describe the multiplicative version, leaving the case of 
an elliptic curve for another occasion. Two sets of free parameters 
enter here: the system of linear differential equations depends on a 
set of \emph{coupling parameters} and their common eigen space 
decomposition is indexed by \emph{spectral parameters}.
We focus on the situation when there is a subsystem of the 
hypergeometric system that is associated to the 
\emph{reflection representation} of the affine Hecke algebra;
this amounts to a particular choice of the spectral parameters. 
This turns out to define a projective structure on the space of 
$W$-orbits  of the toric arrangement complement in question.
For a certain range of the coupling parameters this even yields a complex 
hyperbolic structure on that orbit space. 
It is usually incomplete, but we find that when, in the terminology of 
\cite{CHL03}, the coupling parameters satisfy the \emph{Schwarz conditions},
the situation is as nice as one could hope for: we thus obtain a handful 
of cases for which the hyperbolic structure can be 
completed to a quotient of a complex ball by a discrete group of 
automorphisms acting with 
cofinite invariant volume. We sketch these results in Section 3 with the intention
of treating this material more thoroughly elsewhere.

We found it worthwhile to discuss in some detail (in the final Section 4) 
the example where $R$ is of type $E_8$ 
and the coupling parameter equals $\frac{1}{6}$, because its 
geometry is linked to the moduli spaces of rational elliptic surfaces and 
degree 12 divisors on the projective line \cite{HL02}.
Similar stories can be told for $E_6$ with coupling parameter $\frac{1}{3}$ 
and $E_7$ with coupling parameter $\frac{1}{4}$ in relation to moduli 
spaces of cubic surfaces and quartic curves respectively \cite{ACT02}, \cite{K00}. But
despite these beautiful examples, the picture  of the geometry behind the 
general Calogero-Moser system (with general $R$ and arbitrary spectral and
coupling parameters) remains
unclear. Recent papers on the classical Calogero-Moser system point in the
direction of the classical Hitchin system \cite{HM99}, \cite{BN03}. 
A natural continuation along these lines seems to lead to a further 
analysis of the quantization of Hitchin's fibration \cite{BD96}, \cite{BD97}.

For the reader's convenience we included a brief review of the theory of 
trigonometric Dunkl operators 
in Section 2. 

\medskip\noindent
We dedicate this paper to Gerrit van Dijk on the occasion of his 
65$^{th}$ birthday. One of us (GH) was the first PhD student of Gerrit van 
Dijk at Leiden University, and as such he is grateful for the excellent 
mathematical training he received from Gerrit.

\section{Quantization of the trigonometric Calogero-Moser system}
Let $\mathfrak{a}$ be a real vector space of dimension $n$, 
and let $R \subset\mathfrak{a}^*$ be a possibly nonreduced root 
system spanning $\mathfrak{a}^*$. 
Let $R^\vee \subset \mathfrak{a}$ be the dual root system and write 
$Q^\vee = \Z R^\vee \subset \mathfrak{a}$ for the
coroot lattice of $R$ so that $P := \Hom (Q^\vee,\Z) \subset 
\mathfrak{a}^*$ is the weight lattice of $R$. 
Then $H = \Hom(P,\C^\times)$ is a complex
torus whose character lattice $\Hom (H,\C^\times)$ can be identified with
$P$ and  whose Lie algebra $\Lie (H)=\mathfrak{h}$
can be identified with  $\mathfrak{a}_\C$. Notice that the algebra 
$\C [H]$ of regular functions on $H$ can be identified with the group 
algebra of $P$,
because it has $\{ e^\mu\}_{\mu \in P}$ as an additive basis and the 
multiplication law is simply $e^\mu e^\nu = e^{\mu+\nu}$.

For $\alpha\in R$, we denote by $s_\alpha : 
\lambda \mapsto \lambda - \lambda (\alpha^\vee)\alpha$ 
the corresponding reflection of $\mathfrak{h}^*$. These reflections 
generated the Weyl group $W\subset \GL (\mathfrak{h}^*)$ of $R$.
By duality $W$ also acts on $\mathfrak{h}$ and $H$. 
The reflection $s_\alpha$ fixes the hypertorus 
$H_\alpha = \{ h \in H\, :\, e^\alpha (h) = 1 \}$ pointwise. 
We put $H^\circ= H - \bigcup_{\alpha\in R} H_\alpha$. It is clear that 
the algebra 
$\C[H^\circ]$ of regular functions on $H^\circ$ is obtained from $\C [H]$ 
by inverting the elements $1-e^{-\alpha}$,  $\alpha \in R$. 
For an element $p$ of the symmetric algebra of 
$\mathfrak{h}$, $\sym (\mathfrak{h})$, we denote by
$\partial (p)$ the corresponding translation invariant differential
operator on $H$. So if we identify $\sym(\mathfrak{h})$ with the algebra of polynomial
functions on $\mathfrak{h}^*$, $\C[\mathfrak{h}^*]$, then $\partial (p) e^\mu = p(\mu) e^\mu$
for $\mu \in P$. The algebra of differential operators on $H^\circ$ can (as a 
$\C[H^\circ]$-module) be identified with $\C[H^\circ]\otimes_\C\sym (\mathfrak{h})$; 
we shall denote it by $\D [H^\circ]$. Now $\C [H^\circ]$ is a natural left module for
$\D[H^\circ]$. The group algebra $\C[W]$ of $W$ also acts on $\C[H^\circ]$ 
on the left and together they turn $\C[H^\circ]$ into a left module
for the semidirect tensor product $\D [H^\circ] \otimes_\C \C[W]$ with
multiplication law 
\[
D\otimes w \cdot D'\otimes w' := (DwD'w^{-1}) \otimes ww'.
\]
It is easy to see that $\D[H^\circ] \otimes \C[W]$ acts faithfully on 
$\C [H^\circ]$.

Choose a nonzero $W$-invariant inner product $(\, ,\,)$ on $\mathfrak{a}^*$ 
and denote by $C\in\sym^2\mathfrak{h}$ the associated (complexified) 
quadratic form. Then $\partial (C)$ is the Laplace operator on 
$H$ characterized by $\partial (C)(e^\mu)= C(\mu)e^\mu$, $\mu\in P$.

\begin{definition}\label{def:qhamilton}
Let $\kbold: \alpha\in R\mapsto k_\alpha\in\C$ be a $W$-invariant map, 
extended by zero to all of $\mathfrak{a}$. 
Then the \emph{quantum Hamiltonian} of the trigonometric Calogero-Moser system
with coupling vector $\kbold\in \C^R$ is the element of $\D [H^\circ]$ given by
\begin{equation*}
H_\kbold = \partial (C) + \half\sum_{\alpha\in R} \frac{k_\alpha(1-k_\alpha
- 2k_{2\alpha})(\alpha, \alpha)}
{(e^{\alpha /2} - e^{-\alpha/2})^2}.
\end{equation*}
\end{definition}

If the coupling vector is real, then the quantum Hamiltonian becomes on
the compact subtorus $T:=\Hom (P, U(1))$  of $H$ the
Schr\"odinger operator for a single particle on $T$ repelled or 
attracted along the subtori $T_\alpha = T \cap H_\alpha$ by an inverse
square distance potential.

Let us fix a set of positive roots $R_+$ in $R$.
The complete integrability of the Calogero-Moser system it best understood 
after formal conjugation by
\begin{equation*}
\delta_\kbold^{\frac{1}{2}} = \prod_{\alpha > 0} \left ( e^{\alpha/2} -
e^{-\alpha/2} \right )^{k_\alpha}
= e^{\rho_\kbold}\prod_{\alpha > 0} (1-e^{-\alpha})^{k_\alpha},
\end{equation*}
where $\rho_\kbold = \frac{1}{2}\sum_{\alpha > 0} k_\alpha \alpha$. We get 
\begin{equation*}
\delta _\kbold^{-\frac{1}{2}}H_\kbold \delta _\kbold^{\frac{1}{2}} = 
L_\kbold + (\rho_\kbold,\rho_\kbold),
\end{equation*}
where 
\begin{equation*}
L_\kbold = \partial (C) + \half \sum_{\alpha > 0} k_\alpha (\alpha, \alpha)
\frac{1+e^{-\alpha}}{1 - e^{-\alpha}}\partial
(\alpha^\vee).
\end{equation*}

Differential operators in $\D [H^\circ]$ commuting with $L_\kbold$ yield after
conjugation by $\delta _\kbold^{\frac{1}{2}}$ differential operators in 
$\D[H^\circ]$ commuting with $H_\kbold$ and vice versa. The advantage of 
$L_\kbold$ over $H_\kbold$ will become clear later.

\begin{definition}\label{def_2.2}
For $\xi \in \mathfrak{h}$ the {\em trigonometric Dunkl operator} 
$T_\kbold(\xi )\in\D[H^\circ] \otimes_\C \C[W]$ with coupling vector $\kbold$ is 
\begin{equation*}
T_\kbold(\xi) =(\partial (\xi)-\rho_\kbold (\xi)) \otimes 1 + \sum_{\alpha >0}\frac{k_\alpha
\alpha (\xi)}{1-e^{-\alpha}} \otimes (1-s_\alpha).
\end{equation*}
\end{definition}
This is a first order operator which preserves the subalgebra
$\C[H]$ of $\C [H^\circ]$, because the  divided difference operator
$(1-e^{-\alpha})^{-1} \otimes (1-s_\alpha)$ has that property.

Let us denote by a bar the antilinear involution of $\C[H]$ which sends 
$e^{\mu}$ to $e^{-\mu}$ ($\mu \in P$).
If the coupling vector $\kbold$ lies in $\N^R$, then 
$\delta _\kbold^{\frac{1}{2}} \overline{\delta} _\kbold^{\frac{1}{2}})$ 
is a well-defined element of $\C [H]$ and we have a hermitian inner 
product on $\C[H]$ defined by 
\begin{equation*}
(f,g)\in \C[H]^2\mapsto (f,g)_\kbold:= \text{constant term of }
f\, \overline g \, \delta _\kbold^{\frac{1}{2}} \overline{\delta}_\kbold^{\frac{1}{2}}) /|W|.
\end{equation*} 
It is easy to check that the Dunkl operators are hermitian with respect to 
$(\, ,\, )_\kbold$:
\begin{equation*}
\left( T_\kbold(\xi)f,g \right)_\kbold = 
\left( f, \overline{T_\kbold(\xi)}g \right)_\kbold
=\left ( f, T_\kbold(\overline \xi)g \right)_\kbold,
\end{equation*} 
where the last bar denotes complex conjugation on $\mathfrak{h}$ with 
respect to the real form $\mathfrak{a}$.

Recall that the usual partial ordering $\le$ on $\mathfrak{h}^*$ is defined by $\mu \le \nu$ if
and only if $\nu-\mu \in \N R_+$. Define a new partial ordering $\le_+$ on $P$
as follows. Let $P_+$ be the cone of dominant weights. 
If for $\mu \in P$, $\mu_+$ denotes the unique dominant weight 
in the orbit $W\mu$, then we let 
\begin{equation*}
\mu \le_+ \nu \mbox{ if either } \mu_+ < \nu_+ \mbox{ or } \mu_+ = \nu_+ \mbox{
and } \nu \le \mu.
\end{equation*}
So the smallest element of the orbit $W\mu$ is $\mu_+$ and if $w_0$ 
is the longest element in $W$, then  $w_0 \mu_+$ is the largest element. 
It is now easy to check that Dunkl operators are triangular on $\C [H]$ 
with respect to the basis $e^\mu$ partially ordered by $\le_+$.
The Gram-Schmidt process applied to the monomial basis of $\C [H]$ yields 
a new basis $\{ E_\kbold(\mu)\}_{\mu \in P}$ of $\C[H]$,  characterized by 
\begin{align*}
E_\kbold(\mu) &= e^\mu + \text{lower order terms relative to $\le_+$},\\ 
\left ( E_\kbold(\mu),e^\nu \right )&= 0\quad \text{for all
$\nu\in P$ with $\nu <_+ \mu$}.
\end{align*}
Clearly, the Dunkl operators are also triangular with respect to the basis 
$E_\kbold(\mu)$ partially ordered by $\le_+$. Let $\varepsilon : \R \to \{ \pm 1 \}$ be
given by $\varepsilon (x) = 1$ if $x > 0$ and $\varepsilon (x) = -1$ if $x \le
0$. For $\mu \in P$ let $\widetilde \mu \in \mathfrak{h}^*$ be given by
\begin{equation}\label{eq_2.12}
\widetilde\mu = \mu + \half\sum_{\alpha > 0} k_\alpha \varepsilon \left (
\mu (\alpha^\vee) \right ) \alpha.
\end{equation}
For example for $\mu \in P_+$ dominant and regular $\widetilde\mu = \mu +
\rho_\kbold$, whereas for $\mu \in P_-$ antidominant $\widetilde\mu =\mu -
\rho _\kbold$.
Then a direct computation gives for all $\mu \in P$
\begin{equation*}
T_\kbold(\xi)e^\mu = \widetilde\mu (\xi )e^\mu + \text{lower order terms} ,
\end{equation*}
which by the characterizing property of the basis 
$\{ E_\kbold(\mu)\}_{\mu \in P}$ and the triangularity of $T_\kbold(\xi)$ yields
\begin{equation*}
T_\kbold(\xi) E_\kbold(\mu) = \widetilde \mu (\xi) E_\kbold(\mu) + 
\text{lower order terms.}
\end{equation*}
Using the hermiticity of Dunkl operators and the
characterizing property of the basis $\{ E_\kbold(\mu)\}_{\mu \in P}$ the lower order 
terms are easily found to vanish and we get
\begin{equation}\label{eq_2.14}
T_\kbold(\xi) E_\kbold(\mu) = \widetilde \mu (\xi) E_\kbold(\mu)\quad \text{for all $\mu
\in P$}.
\end{equation}
Indeed triangular plus symmetric becomes diagonal. We conclude:

\begin{theorem}
For fixed vector $\kbold$ of coupling constants the Dunkl operators 
$T_\kbold(\xi)$ for $\xi\in \mathfrak{h}$ mutually commute. In other words,
$T_\kbold$ extends to 
an algebra homomorphism $T_\kbold: \sym(\mathfrak{h})\to \D[H^\circ]\otimes \C[W]$.
\end{theorem}

We denote the coefficient of $w\in W$ in $T_\kbold(p)$ by $D_\kbold(p,w)\in\D (H^\circ)$, 
so that 
\[
T_\kbold(p) = \sum_{w\in W} D_\kbold(w,p)\otimes w.
\] 
Our choice of positive roots determines (and is determined by) a
set of simple roots $\alpha_1 ,\ldots , \alpha_n$. We shall write
$s_i$ for $s_{\alpha_i}$. A rank one calculation gives for $i=1,\dots, n$
\begin{equation*}
s_i T_\kbold(\xi) - T_\kbold(s_i\xi ) s_i = - 
(k_{\alpha_i} + 2k_{2\alpha_i})\alpha_i (\xi)\;\text{ for all 
$\xi\in\mathfrak{h}$}
\end{equation*}
(recall that we put $k_\beta = 0$ if $\beta \notin R$). 
With induction on the degree this yields
\begin{equation*}
s_i T_\kbold(p) - T_\kbold(s_i p) s_i = -(k_{\alpha_i} + 2k_{2\alpha_i})
T_\kbold \left (\frac{p-s_ip}{\alpha_i^\vee}\right ) 
\;\text{ for all $p\in\sym(\mathfrak{h})$.}
\end{equation*}
In particular we find that for 
$q\in \sym(\mathfrak{h})^W$, $T_\kbold(q)$ commutes with $W$. So if we then write
$T_\kbold(q) = \sum_{w\in W} D_\kbold(w,q) \otimes w$, then
\begin{equation}\label{eq_2.141}
D_\kbold(q) = \sum_{w\in W} D_\kbold(w,q) \in \D [H^\circ]^W
\end{equation}
is the unique $W$-invariant differential operator which has the same restriction
to $\C[H]^W$ as the Dunkl operator $T_\kbold(q)$. In particular, $D_\kbold(q)$ preserves
$\C[H]^W$ because $T_\kbold(q)$ does. Since the Dunkl operators commute, we find:

\begin{theorem}\label{theo_2.3}
The map $D_\kbold:\sym(\mathfrak{h})^W\to \D[H^\circ]^W$
is a $\C$-algebra homomorphism. We have $D_\kbold(C)=L_\kbold + (\rho_\kbold,\rho_\kbold)$
and the image of $D_\kbold$ preserves $\C[H]^W$.  
\end{theorem}

Since the image of $D_\kbold$ contains the
operator $L_\kbold$, conjugation by $\delta _\kbold^{\frac{1}{2}}$ of this algebra
immediately gives the quantum complete integrability of the
trigonometric Calogero-Moser system. The classical complete integrability
follows as in the thesis of Opdam. Indeed one puts  
$g_\alpha^2 = k_\alpha (1-k_\alpha - 2 k_{2\alpha})$ and observes 
that the differential operators 
$\delta_\kbold^{\frac{1}{2}} 
\circ D_\kbold(p) \circ\delta_\kbold^{-\frac{1}{2}}$ for 
$p\in\sym(\mathfrak{h})^W$
depend in a polynomial way on the $g_\alpha$ for $\alpha \in R$. 
The grading on $\sym(\mathfrak{h})$ is extended to the coupling constants 
$g = (g_\alpha)$ by giving each of these degree one. A symbol calculus for this total degree, which
amounts to the classical limit, yields the classical complete integrability.

\medskip\noindent
It was observed by Matsuo and by Cherednik that the eigenvalue problem for the
Dunkl operators can be rewritten as an integrable connection of
Knizhnik-Zamolodchikov type. For an exposition of this result as well as further
references of the theory explained in this section we refer to the survey
articles of Heckman \cite{H97} and of Opdam \cite{O00}.

\section{A special eigenstate of the trigonometric Calogero-Moser system}
\emph{Here and in the rest of the paper we assume that $R$ is a 
reduced root system.} We write $k_i = k_{\alpha_i}$ for $i = 1,\ldots , n$.
Note that the vector $\rho_\kbold$ is now characterized by 
$\rho_\kbold(\alpha_i^\vee)=k_i$.

The simultaneous eigenvalue problem for the trigonometric Calogero-Moser system
is obtained from the system
\begin{equation}\label{eq_3.1}
D_\kbold(p)F = p(\lambda)F \; \mbox{ for all } p \in \sym(\mathfrak{h})^W
\end{equation}
by conjugation with the weight function $\delta _\kbold^{\frac{1}{2}}$. Here
$\lambda \in \mathfrak{h}^*$ is fixed, and called the \emph{spectral parameter} of the
system (\ref{eq_3.1}). For the (admittedly, nonreduced) root system 
$BC_1$, (\ref{eq_3.1}) boils down to the Euler-Gauss hypergeometric 
equation, and therefore (\ref{eq_3.1}) is called the 
\emph{hypergeometric system} associated to the
root system $R$. The above system has regular singular points at infinity on
$H^\circ$. The dimension of the local solution space equals $|W|$ and the
monodromy representation of the fundamental group of $W\backslash H^\circ$,
which is the affine Artin group associated with $R$, factors through the
affine Hecke algebra with quadratic relations
\begin{equation}\label{eq_3.2}
\left (M_i -1\right ) \left ( M_i + e^{-2\pi i k_i} \right ) = 0.
\end{equation}
Here $M_i$ is the monodromy operator associated to the analytic 
continuation of a solution starting in the image of the 
positive chamber under the exponential map and 
making a half turn around the subtorus $H_i$ corresponding to the 
simple root $\alpha_i$. 

\medskip
Following Harish-Chandra we substitute into the hypergeometric system 
(\ref{eq_3.1}) a formal series of the form
\begin{equation*}
F = \sum_{\nu \le \mu} a_\nu e^\nu \quad\text{with}\quad a_\mu = 1.
\end{equation*}
With the help of (\ref{eq_2.12}) and (\ref{eq_2.14}) we see that the 
$exponents \mu \in\mathfrak{h}^*$ for which such formal solutions exist 
must satisfy the \emph{indicial equation} 
$p(\mu + \rho_\kbold) = p (\lambda)$ for all $p \in S \mathfrak{h}^W$. 
This amounts to
\begin{equation*}
\mu \in W \lambda - \rho_\kbold.
\end{equation*}
Choosing the base point in the positive chamber the eigenvalues 
of the monodromy operators of the fundamental group of $H^\circ$
are easily read of from these exponents.
Now, if the spectral parameter $\lambda \in \mathfrak{h}^*$ is such that all
points $w \lambda$ for $w \in W$ are distinct modulo the root lattice $Q$, 
i.e., if $\lambda (\alpha^\vee) \notin \Z$ for all $\alpha \in R$,
then the hypergeometric system (\ref{eq_3.1}) has a formal solution basis. 
Under this nonresonance condition it is easy to see that
the monodromy representation is reducible if and only if
\begin{equation}\label{eq_3.7}
\lambda (\alpha^\vee) + k_\alpha \in \Z \; \mbox{ for some } \alpha \in R.
\end{equation}
In fact this criterion remains valid for $\lambda$ resonant. We shall stick to
the notation
\begin{equation}\label{eq_3.8}
\lambda = \mu + \rho _\kbold.
\end{equation}
So $\mu$ is an exponent at infinity (in the direction of the positive chamber), whereas $\lambda$,
obtained by a shift over $\rho_\kbold$, is the spectral parameter of our
hypergeometric system.

If $\mu \in P_+$ and $k_\alpha >0$ for all $\alpha \in R$, then the hypergeometric
system (\ref{eq_3.1}) has solutions in $\C[H]^W$ with leading exponent $\mu =
\lambda - \rho_\kbold$. This is an extreme case of reducibility since (remember that
$\rho_\kbold(\alpha_i^\vee) = k_i$) we have $\lambda (-\alpha_i^\vee) + k_i = \mu
(-\alpha_i^\vee) \in -\N$ for $i = 1,\ldots ,n$. Hence for generic $\kbold$ the
exponents $\mu \in P_+$ are intersection points of normal crossing reducibility
hyperplanes (\ref{eq_3.7}). Clearly the solutions in $\C[H]^W$ are fixed under
the monodromy representation.

The goal of this section is to discuss the next best case for which the
hypergeometric system admits a monodromy invariant subspace of solutions of
dimension $n+1$ that affords the reflection representation of the affine Artin
group as its monodromy. These special solutions of (\ref{eq_3.1}) 
satisfy a bigger set of differential equations, a set which we baptize
the \emph{special hypergeometric system}. It
turns out to be generated by quadratic equations and to give rise to a
projective structure on $H^\circ$.

\medskip\noindent
\emph{From now on we assume that the reduced root system $R$ is irreducible.} 
We denote by $\varpi_1,\ldots , \varpi_n$ the fundamental weights in $P_+$. In dealing with
particular cases the numbering will be the same as in the Bourbaki text on root
systems \cite{B68}. This means in particular that $R=W\{\alpha_1,\alpha_n\}$. 

\medskip\noindent
We shall describe for each root system a sequence 
$(\mu_1,\ldots , \mu_{n+1})$ in $\mathfrak{h}^*$ given up to order
(repetition is not excluded) which describes
the exponents of the special hypergeometric system with multiplicities. 
This will be done in a case by case manner. 
The points $\mu_i=\mu_{i,\kbold}$ will depend linearly on $\kbold$ and 
in the case of type $A_n$ on an additional parameter as well. 
Since $\kbold$ takes on $R$ at most two values, 
namely $k_1$ and $k_n$, we find it
convenient to adapt our notation accordingly: we write $k$ for $k_1$ and 
if $\alpha_n\notin W\alpha_1)$ (i.e., $R$ of type $B$, $C$, $F$ or $G$), 
$k'$ stands for $k_n$. In the $A_n$-case, $k'$ will have a different 
meaning, but we will then abuse notation, by letting $\kbold$
then stand for the pair $(k=k_1,k')$.
We shall first discuss a root system whose Dynkin diagram is a chain, 
i.e., one of type $A$, $B$, $C$, $F$ or $G$, and
subsequently deal with those of type $D$ or $E$. 

For $R$ of type $A_n$ $(n \ge 2)$, the extended Dynkin diagram has a 
loop causing an additional modulus in the affine reflection representations.
This accounts for an additional parameter $k' \in \C$. We put 
\begin{equation}\label{eq_3.9}
\mu_i = (x-ik)\varpi_{i-1} + ((i-1)k-x)\varpi_i
\end{equation}
for $i = 1,\ldots , n+1$ with the convention $\varpi_0 = \varpi_{n+1} = 0$. Here
we take $x,y \in \C$ related to $k,k' \in \C$ by
\begin{equation*}
x+y = (n+1)k, \; x-y = (n+1)k'
\end{equation*}
so that $\mu_1 = -x \varpi_1$, $\mu_{n+1} = -y \varpi_n$. Since $\alpha_i =
-\varpi_{i-1} + 2 \varpi_i - \varpi_{i+1}$ it is clear that
\begin{equation*}
\mu_{i+1} - \mu_i = (x-ik)\alpha_i \; \mbox{ for } i = 1,\ldots ,n .
\end{equation*}
Keeping in mind (\ref{eq_3.8}) it is clear from (\ref{eq_3.9}) that $\lambda_i
(\alpha_i^\vee) = ik-x$ and so the last equality yields
\begin{equation*}
\lambda_{i+1} = s_i (\lambda_i) \; \mbox{ for } i = 1,\ldots ,n.
\end{equation*}

For $R$ of type $B_n$ we still take (\ref{eq_3.9}) for $i = 1,\ldots, n-1$ 
and in accordance with (\ref{eq_3.8}) and the above identity, we put
\[
\begin{split}
\mu_n &= \mu_{n-1} - ((n-1)k-x)\alpha_{n-1} \\
&= (x-nk)\varpi_{n-1} + (2(n-1)k-2x)\varpi_n,\\
\mu_{n+1} &= \mu_n {-} (2(n{-}1)k{+}k'{-}2x)\alpha_n\\ 
&= ((n{-}2)k{+}k'{-}x)\varpi_{n-1} {+}(2x{-}2(n{-}1)k){-}2k')\varpi_n,
\end{split}
\]
with the convention (for $R$ of type $B$, $C$ $F$ or $G$) that $k' = k_n$.
Hence if we take $x = (n-2)k+k'$, $y = 2k$ then we find 
$\mu_1 = -x\varpi_1$ and $\mu_{n+1} = -y \varpi_n$.

Likewise for $R$ of type $C_n$ we find
\begin{align*}
\mu_n &= \mu_{n-1} {-} ((n{-}1)k{-}x)\alpha_{n-1} \\
&= (x{-}nk)\varpi_{n-1} {+}
((n{-}1)k{-}x)\varpi_n,\\
\mu_{n+1} &= \mu_n {-} ((n{-}1)k{+}k'{-}x)\alpha_n\\ 
&=((n{-}2)k {+} 2k'{-}x)\varpi_{n-1} + (x{-}(n{-}1)k {-}2k') \varpi_n,
\end{align*}
and with $x = (n-2)k + 2k'$, $y = k$ we get $\mu_1 = -x\varpi_1$ and $\mu_{n+1} =
-y\varpi_n$. 

In case $R$ is of type $F_4$ we take
\begin{align*}
\mu_1 &= -2x\alpha_1 -3x \alpha_2 -4x \alpha_3 -2x \alpha_4 =
-(k+k')\varpi_1,\\
\mu_2 &= -y\alpha_1 -3x \alpha_2 -4x \alpha_3 -2x\alpha_4 = (-k+k')\varpi_1 -
k'\varpi_2,\\
\mu_3 &= -y \alpha_1 -2y \alpha_2 -4x \alpha_3 -2x \alpha_4 = (-2k+k')\varpi_2 +
(2k-2k')\varpi_3,\\
\mu_4 &= -y \alpha_1 - 2y\alpha_2 - 3y \alpha_3 -2x \alpha_4 = -2k\varpi_3 +
(2k-k')\varpi_4,\\
\mu_5 &= -y \alpha_1 - 2y \alpha_2 - 3y \alpha_3 - y \alpha_4 = -(2k+k')\varpi_4,
\end{align*}
with $x = k+k'$, $y = 2k+k'$. Note that $\mu_2 - \mu_1 = k'\alpha_1$,
$\mu_3-\mu_2 = (-k+k')\alpha_2$, $\mu_4-\mu_3 = (-2k+k')\alpha_3$ and
$\mu_5-\mu_4 = -2k\alpha_4$.

Finally for $R$ of type $G_2$, let
\begin{align*}
\mu_1 &= -2x \alpha_1 - x \alpha_2 = -x \varpi_1,\\
\mu_2 &= -3y \alpha_1 - x \alpha_2 = (x-2k)\varpi_1 + (k-x)\varpi_2,\\
\mu_3 &= -3y \alpha_1 - 2y\alpha_2 = -y\varpi_2,
\end{align*}
with $x = \frac{1}{2} (k+3k')$, $y = \frac{1}{2} (k+k')$. Note that $\mu_2 -
\mu_1 = \frac{1}{2}(-k+k')\alpha_2$. This completes our case by case 
descriptions of the set of special exponents for 
$R$ of type $A$, $B$, $C$, $F$ or $G$.

For $R$ of type $D_n$ we take for $i = 1,\ldots , n-2$ 
\begin{equation*}
\mu_i = (n-2-i) k\varpi_{i-1} - (n-1-i)k\varpi_i, \; \mu_{n-1} =
-2k\varpi_{n-1}, \; \mu_n = -2k\varpi_n,
\end{equation*}
where $\varpi_0 = 0$, as before. With  $x =(n-2)k$ and 
$y = 2k$ we have $\mu_1 =
-x\varpi_1$, $\mu_{n-1} = -y\varpi_{n-1}$, $\mu_n = -y \varpi_n$.

For $R$ of type $E_8$ let 
\begin{align*}
\mu_1 = -3k\varpi_1, & \mu_5 = -2k\varpi_5 + k\varpi_6,\\
\mu_2 = -2k \varpi_2, & \mu_6 = -3k\varpi_6 + 2k\varpi_7,\\
\mu_3 = k\varpi_1 - 2k\varpi_3, & \mu_7 = -4k \varpi_7 + 3k\varpi_8,\\
\mu_4 = -k\varpi_4, & \mu_8 = -5k\varpi_8.
\end{align*}
The case $R$ of type $E_7$ is obtained by forgetting $\mu_8$ and putting
$\varpi_8 = 0$ in the expression for $\mu_7$ and likewise for $E_6$.
Note that for both $D_n$ and $E_n$
\[
\lambda (\alpha_i^\vee) = - d(i)k
\]
with $d(i)$ the distance in the Dynkin diagram from the i$^{th}$ node to the
triple node. For type $D_n$ and $E_n$ put $\mu_{n+1} = \mu_{n-2}$ and $\mu_{n+1}
= \mu_4$ respectively, so $\mu_i$ occurs with multiplicity two if $i$ 
labels the triple node.  
All in all, we have now given for each irreducible root system
$R$ of rank $n$ a set $\{ \mu_1 ,\ldots , \mu_{n+1} \}$ in $\mathfrak{h}^*$
depending linearly on $k$ resp.\ on $(k,k')$ depending on whether or not
the type is $D$ or $E$. We shall write $K$ for the space of coupling parameters, so 
$k\in K =\C$ or $(k,k')\in K = \C^2$ respectively.

\begin{definition}
The \emph{special exponents} of $R$ are the terms of  
$(\mu_1 ,\ldots , \mu_{n+1})$ in $\mathfrak{h}^*$, 
counted with multiplicity. 
\end{definition}

Notice that each $\mu_i$ depends linearly on $\kbold \in K$.
We denote the inner product on $\mathfrak{a}$ dual to $(\, ,\, )$ by
$(\, ,\, )^\vee$ and write $C^\vee$ for  the corresponding element of  
$\sym^2(\mathfrak{h}^*)$. 

For the $A_n$ root system we define for every
positive root $\alpha\in\mathfrak{a}^*$ an element $\alpha'\in\mathfrak{a}$ 
as the vector in the characterized by the property that 
$\beta (\alpha ')\ge 0$ for all $\beta\in R_+$ with equality when 
$\beta=\alpha$ or $(\beta ,\alpha )=0$ and with 
$(\alpha',\alpha ')^\vee(\alpha,\alpha )=4\frac{n-1}{n+1}$
In terms of the standard description in Bourbaki \cite{B68}:  
if $\alpha^\vee= e_i-e_j$ with $i<j$, then $\alpha'$ is the image of
$e_i + e_j$ in $\R^{n+1}/\text{main diagonal}$.

\begin{proposition}\label{prop_3.2}
The special exponents are exactly the solutions in $\mu$ 
of the following quadratic equation in $\sym^2(\mathfrak{h}^*)$: 
\begin{align*}
\mu^2& + \tfrac{1}{2}\sum_{\alpha>0}\mu(k\alpha^\vee+k'\alpha')\alpha^2 + 
a(k,k')C^\vee = 0\quad  \text{for type $A_n$, $n \ge 2$,}\\
\mu^2& + \tfrac{1}{2}\sum_{\alpha>0}\mu (k_\alpha\alpha^\vee)
\alpha^2+a(k,k')C^\vee=0
\quad \text{for type $B$, $C$, $F$ or $G$,}\\
\mu^2& + \tfrac{1}{2}\sum_{\alpha>0}\mu (k_\alpha \alpha^\vee)\alpha^2+
a(k)C^\vee=0\quad \text{for type $D$ or $E$}.
\end{align*}
The function $a$ is the quadratic polynomial in the coupling constants:
for $R$ of type $A$, $B$, $C$, $F$ or $G$ it is given by
$a(k,k') = xy(\varpi_1,\varpi_n) = (\mu_1, \mu_{n+1})$ and for type 
$D_n$, $E_6$, $E_7$, $E_8$ (with inner product normalized by 
$(\alpha, \alpha) = 2$) by $a(k) = (n-2) k^2$ and $a(k) = 6k^2, 12k^2, 30k^2$
respectively.
\end{proposition}

For $\kbold\in\ K$, we define a linear map
\begin{equation}\label{eq_3.141}
D_\kbold: \sym^2(\mathfrak{h})\to \D [H^\circ]
\end{equation}
by letting its value on the monomial $\xi\eta$ (with $\xi,\eta\in\mathfrak{h}$) be
\[
\partial(\xi)\partial(\eta) +  \tfrac{1}{2} \sum_{\alpha > 0}
\alpha(\xi)\alpha(\eta) \frac{1+e^{-\alpha}}{1-e^{-\alpha}} k_\alpha\partial
(\alpha^\vee) + \rho_\kbold(\xi)\rho_\kbold(\eta) 
\]
if the root system is not of type $A_n$ and 
\[
\partial(\xi)\partial(\eta) +  \tfrac{1}{2} \sum_{\alpha > 0}
\alpha(\xi)\alpha(\eta) \left(\frac{1+e^{-\alpha}}{1-e^{-\alpha}}
k\partial(\alpha^\vee)+k'\partial(\alpha')\right) + \rho_\kbold(\xi)\rho_\kbold(\eta) 
\]
if the root system is of type $A_n$, $n\ge 2$. 
Observe that the notations in (\ref{eq_2.141}) and (\ref{eq_3.141})
are compatible for $p\in \sym^2(\mathfrak{h})^W$.

\begin{theorem}\label{theo_3.3}
The system of differential equations
\begin{equation}\label{eq_3.27}
D_\kbold(p) F = (p(\rho_\kbold)-a(\kbold)(C^\vee,p))F \quad \text{for all $p\in \sym^2(\mathfrak{h})$},
\end{equation}
enjoys the following properties:
\begin{enumerate}
\item[(i)] the system is integrable on $H^\circ$ and invariant under $W$, 
\item[(ii)] solutions of (\ref{eq_3.27}) are also solutions of the
general hypergeometric system (\ref{eq_3.1}) with spectral parameter 
$\lambda_i=\mu_i(\kbold)+\rho_\kbold$, where $\mu_i$ is a special exponent,
\item[(iii)] the dimension of the local solution space of (\ref{eq_3.27}) 
is equal to $n+1$, and the monodromy of the associated rank $n+1$ local 
system over $W\backslash H^\circ$ factors through 
the reflection representation of the affine Hecke algebra (viewed as a
quotient of the group algebra of the fundamental group of 
$W\backslash H^\circ$): $M_i$ has eigenvalue 1 with
multiplicity $n$ and eigenvalue $-e^{2\pi ik_i}$ with multiplicity 1, in
accordance with (\ref{eq_3.2}).
\end{enumerate}
\end{theorem}

We call (\ref{eq_3.27}) the {\em special hypergeometric system} with
parameter $\kbold$.

\begin{theorem}\label{theo_3.4}
For real $\kbold \in K$ the special hypergeometric system (\ref{eq_3.27}) admits a
monodromy invariant hermitian form $h_k$. This is nondegenerate and has 
Lorentzian signature if $\kbold\in K_+$, where
\[
K_+ := 
\begin{cases}
(0,\frac{1}{n-2})&\text{ for type $D_n$,}\\ 
(0,\frac{1}{n-3})&\text{ for type $E_n$,}\\
\{ (k,k') \in (-\half,\half)^2\, :\, 
0 < x < 1\, ;\, 0 < y < 1 \}&\text{otherwise.}\\ 
\end{cases}
\]
\end{theorem}

\begin{lemma}\label{lem_3.5}
For a basis $F_0,\ldots , F_n$ of a local solution space of (\ref{eq_3.27}) and
a basis $\xi_1, \ldots , \xi_n$ of $\mathfrak{h}$ the Wronskian is 
(up to a constant depending on choices of bases) given by
\begin{equation*}
\det \begin{pmatrix}
F_0 & F_1 & \ldots & F_n\\
\partial(\xi_1)F_0 & \partial(\xi_1)F_1 & \ldots & \partial(\xi_1)F_n \\
\vdots & \vdots && \vdots\\
\partial(\xi_n)F_0 & \partial(\xi_n)F_1 & \ldots & \partial (\xi_n)F_n 
\end{pmatrix}
= \delta _\kbold^{-1}.
\end{equation*}
\end{lemma}

Fix a base point $z_0 \in W\backslash H^\circ$, and let 
$\kbold \in K_+$. We denote the
projectivized monodromy group of (\ref{eq_3.27}) by 
$\Gamma = \Gamma_\kbold$ and by
$\widetilde{W\backslash H^\circ} \to W\backslash H^\circ$ the 
corresponding $\Gamma$-covering. It follows from 
Lemma \ref{lem_3.5} that the projectivized evaluation map
given by $[F_0(z) \, :\, \ldots : F_n(z)]$
becomes a local isomorphism
\begin{equation*}
\Pev: \widetilde{W\backslash H^\circ}{\longrightarrow} \mP^n(\C)
\end{equation*}
and so equation (\ref{eq_3.27}) gives a projective structure on $W\backslash
H^\circ$. The monodromy invariant hermitian form yields a complex hyperbolic
ball $\mB^n \subset \mP^n$ on which $\Gamma$ acts by biholomorphic automorphisms.

\begin{definition}\label{def_3.6}
We say that $\kbold \in K_+$ is \emph{admissible} if the image of $\Pev$ is
contained in $\mB^n$, so that we have a map $\Pev:\widetilde{W\backslash H^\circ}\to\mB^n$.
\end{definition}

\noindent Admissibility was obtained for the additive setting
\cite{CHL03} and we expect that its proof can be carried 
over to the multiplicative setting. For the rest of this section we assume that the
parameter $\kbold$ is admissible.

It is our goal to find conditions on $\kbold$ such that $\Gamma$ acts on $\mB^n$ as a
lattice group in $PU(n,1)$, so that the commutative diagram
\begin{equation}\label{eq_3.33}
\begin{array}{ccc}
\widetilde{W\backslash H^\circ} & \overset{\Pev}{\longrightarrow} & \mB^n\\
\downarrow &&\downarrow\\
W\backslash H^\circ &\longrightarrow & \Gamma\backslash \mB^n
\end{array}
\end{equation}
lives in the complex-analytic category.
Notice that in the bottom line the lefthand side is in the realm of 
{\em projective} geometry, whereas the right hand side concerns 
\emph{hyperbolic} geometry.
Even if $\Gamma$ does not act as a lattice group on $\mB^n$ the space
$W\backslash H^\circ$ always acquires a hyperbolic metric under the
admissibility assumption $\kbold \in K_+$. If it happens to be complete, then $\Pev :
\widetilde{W\backslash H^\circ} \to \mB^n$ is a covering map, and hence a
holomorphic isomorphism, since $\mB^n$ is simply connected. So 
$W\backslash H^\circ$ has then the structure of a ball quotient. 

However it is rare that the hyperbolic metric on $W\backslash H^\circ$ is
complete. In the one-dimensional case of the Euler-Gauss hypergeometric 
equation there is
only one case, namely when all exponent differences are zero and $\Pev$ gives the
uniformization of $\mP^1 - \{ 0,1,\infty \}$ by the unit disc.
To remedy the above argument one seeks a (partial) stratified compactification on
the projective side of diagram (\ref{eq_3.33}) and an extension of $\Pev$ over the
associated boundary as a locally biholomorphic map, such that the hyperbolic metric on the
compactified space is complete. Now one can argue as before. To find the correct
compactification on the projective side is a very delicate problem in birational
geometry, and its mere existence imposes strong restrictions on $\kbold \in K_+$,
called the \emph{Schwarz conditions} in the additive case \cite{CHL03}.
A direct approach to find the desired compactification is due to Deligne and
Mostow \cite{DM86}, and uses the  geometric invariant theory of binary forms. 
Its scope is limited to the case of classical root systems, but it does construct on the
projective geometric side the compactification that corresponds to the
Baily-Borel compactification on the ball quotient side.

An alternative approach, followed in \cite{CHL03}, is to work with a 
rather easily described compactification that dominates the searched 
compactification. The map $\Pev$ extends to a holomorphic submersion 
onto $\mB^n$ with connected components of fibers being compact, and 
one can apply Stein factorization to it; the Stein factor then gives 
the desired compactification. Although this procedure is hardly 
constructive in the projective category, the outcome is that 
the complement of the image of 
$\Pev$ in $\mB^n$ is a locally finite union of hyperballs (a so-called Heegner
divisor), and typically the group $\Gamma$ contains reflections with 
the hyperballs as mirrors.

\begin{question}\label{ques_3.7}
Is there a modular interpretation of diagram (\ref{eq_3.33}), such that the
projectivized evaluation map $\Pev$ gets the interpretation of a period map?
\end{question}
For $R$ a classical root system the answer to the question is positive, and
given by the theory of Deligne and Mostow \cite{DM86}. For $R$ of type $E$
one encounters moduli spaces of Del Pezzo surfaces and suitable period maps 
on these \cite{ACT02}, \cite{K00}, \cite{HL02}. But the picture for general 
root systems $R$ remains unclear. In the next section we shall discuss the 
case $R$ of type $E_8$ with $k=\frac{1}{6}$.

\section{An example: the $E_8$ root system with coupling constant $\frac{1}{6}$}
The Schwarz conditions for the classical Euler-Gauss hypergeometric equation
require that at the three singular points $0,1,\infty$ the exponent differences
are the reciprocal of integers. Once the correct stratified compactification of
$W\backslash H^\circ$ adapted to the multivariable special hypergeometric system
(\ref{eq_3.27}) is found the Schwarz conditions boil down to similar conditions
as for the one dimensional case, transverse to the codimension strata of the boundary. The special hypergeometric system
(\ref{eq_3.27}) was given as a system on $H^\circ$ invariant under $W$. However
from Theorem (\ref{theo_3.3}) it is clear that the system lives in fact on the adjoint
torus $H'$ (with character lattice $Q = \Z R$), and is invariant under $W' = W
\rtimes \mbox{Aut}(R_+)$ for $R$ not of type $A_n$. For type $A_n$ we have to
impose the condition $k' = 0 \iff x = y$ so $-w_0\mu_1 = \mu_{n+1}$. In
principle working on the quotient space $W' \backslash H'^{\circ}$  of $W
\backslash H^\circ$ allows more solutions to the Schwarz conditions. However in
certain examples there is still more symmetry which is unclear from the root
system perspective. One  ought to think of this as a {\em hidden symmetry} for
the system (\ref{eq_3.27}).

\begin{example}\label{ex_4.1}
For the root system $R$ of type $A_n$ one has a covering map
\begin{equation*}
W'\backslash H'^{\circ} \longrightarrow {\mathcal D}_{n+3}
\end{equation*}
of degree $\binom{n+3}{2}$ with ${\mathcal D}_{n+3}$ 
the moduli space of reduced effective degree $(n+3)$
divisors on the projective line $\mP$. Indeed $H^\circ = \{ (x_1,\ldots ,
x_{n+1}) \in (\C^\times )^{n+1}$; $x_1.\ldots .x_{n+1} = 1$, $x_i \ne x_j$ for
$1 \le i < j \le {n+1}\}$ and $H'^{\circ} = C_{n+1}\backslash H^\circ$ with
$C_{n+1}$ the group of roots of unity of order $(n+1)$ acting diagonally, while $W = S_{n+1}$ and
$W' = S_{n+1} \rtimes S_2$. Writing $H'^{\circ} = \C^\times \backslash \{ (x_0
= 0,x_1,\ldots , x_{n+1}, x_{n+2} = \infty) \in \mP^{n+3}$; $x_i \ne x_j$ for $0
\le i < j \le n+2\}$ the map $(x_0,\ldots , x_{n+2} ) \mapsto \{ x_0 ,\ldots ,
x_{n+2}\}$ yields the above covering map. The special hypergeometric system is the
pull back of the Lauricella system on ${\mathcal D}_{n+3}$ if $(n+3)k = 2$. Under
this hypothesis, the signature of the invariant hermitian form is 
Lorentzian. The Schwarz conditions amount to $(\frac{1}{2} - k )^{-1} = q \in \N$.
This leaves us with the following possibilities
\begin{equation*}
\begin{tabular}{c|c|c|c|c|c}
$n$&1&2&3&5&9\\
\hline 
$q$&$\infty$ & 10&6&4&3
\end{tabular}
\end{equation*}
which are exactly those cases in the table of Deligne and Mostow with all
Lauricella exponents equal (to $k$). In all cases with $q \in 2 \N$ the Schwarz condition
would already be satisfied on ${\mathcal M}_{0,n+3}$. The quotient by $S_{n+3}$ is
therefore only needed to include the last case $n = 9$, $q=3$ which gives
\begin{equation}\label{eq_4.3}
{\mathcal D}_{12} \hookrightarrow \Gamma_9 \backslash \mB^9
\end{equation}
with $\Gamma_9$ generated by order 3 hyperbolic reflections. Their mirrors are
hyperballs in $\mB^9$ whose quotient by $\Gamma_9$ is exactly the boundary of
${\mathcal D}_{12}$ under the map (\ref{eq_4.3}). For $D \in {\mathcal D}_{12}$
represented by $\{ x_0,\ldots , x_{11} \}$ consider the algebraic curve $C$ with
affine equation $y^6 = \prod_{i=0}^{11} (x-x_i)$. The periods
\begin{equation*}
\int_{x_i}^{x_{i+1}} \; \frac{dx}{y}
\end{equation*}
are Lauricella integrals and solutions of (\ref{eq_3.27}) for $R$ of type $A_9$
with $(k,k') = (\frac{1}{6},0)$.
\end{example}

\medskip\noindent
A rational elliptic surface $S$ in Weierstrass form is a projective surface with
affine equation
\begin{equation}\label{eq_4.5}
S : y^2 = 4x^3 - g_2x - g_3
\end{equation}
with $g_2,g_3 \in \C [u,v]$ binary forms of degree 4 and 6 respectively. 
The surface represents an elliptic fibration over
the projective line $\mP$ with coordinate $t = u/v$. The \emph{discriminant
divisor} $D$ of the fibration is the divisor on that line defined by 
$\Delta:=g_2^3 - 27 g_3^2 \in \C[u,v]$; its support gives the singular fibers
of the fibration, but $D$ gives each such fiber a multiplicity.
Since $\Delta$ is a binary form of degree 12, so $D$ is an effective degree 12 divisor on $\mP$.  
We have a zero of multiplicity $1$ precisely 
when the corresponding singular fiber is a nodal rational curve.
Let ${\mathcal M}$ be the moduli space of rational elliptic
surfaces of the form (\ref{eq_4.5}) with reduced discriminant divisor. There is
a natural map
\begin{equation*}
{\mathcal M} \to {\mathcal D}_{12} , \; S  \mapsto D
\end{equation*}
which turns out to be injective. Since $\dim {\mathcal M} = 8$ and 
$\dim {\mathcal D}_{12} 
= 9$ the space ${\mathcal M}$ lies as a hypersurface in ${\mathcal D}_{12}$. 
The main result of \cite{HL02} is the following theorem.

\begin{theorem}\label{theo_4.2}
We have a commutative diagram
\begin{equation*}
\begin{array}{ccc}
{\mathcal M} & \longrightarrow & \Gamma_8\backslash \mB^8\\
\downarrow && \downarrow\\
{\mathcal D}_{12} & \longrightarrow & \Gamma_9\backslash\mB^9
\end{array}
\end{equation*}
with the bottom arrow the Deligne-Mostow period map (\ref{eq_4.3}).
\end{theorem}
Relative to the natural line bundle on ${\mathcal D}_{12}$, the hypersurface ${\mathcal M}$ in
${\mathcal D}_{12}$ has degree 3762 \cite{V99}, and so from a projective point of view
${\mathcal M}$ lies in a wildly bended manner in ${\mathcal D}_{12}$. On the 
other hand, through hyperbolic eyes the embedding is totally geodesic 
(which can be thought of as linear). To understand the connection of Theorem 
\ref{theo_4.2} and the previous section we recall a result of Looijenga \cite{L98}.

\begin{theorem}\label{theor_4.3}
Let ${\mathcal M}'$ be the moduli space of rational elliptic surfaces of the form
(\ref{eq_4.5}) with reduced discriminant, marked with a distinguished singular
nodal rational fiber. Then we have a commutative diagram
\begin{equation}\label{eq_4.9}
\begin{array}{ccc}
{\mathcal M}' & \overset{\cong}{\longrightarrow} & W\backslash H^\circ\\
&\searrow \qquad \swarrow&\\
& {\mathcal M} &
\end{array}
\end{equation}
with ${\mathcal M}' \to {\mathcal M}$ the natural degree 12 covering map given by forgetting the
distinguished fiber, and $W \backslash H^\circ$ the regular toric orbit space
for $R$ of type $E_8$.
\end{theorem}
The description of the isomorphism (\ref{eq_4.9}) is in geometric terms. There
is a distinguished section in the Weierstrass model $S$ over $\mP$ through $(x :
y : z) = (0 : 1 : 0)$ called the zero section. The global sections of $S$ over
$\mP$ form a lattice by addition in fibers, and with respect to minus the
intersection form on $S$, is called the \emph{Mordell-Weil lattice}. 
It turns out for $S$ in ${\mathcal M}$ that the Mordell-Weil lattice is positive 
definite, even, unimodular and of rank 8, and
hence isomorphic to the $E_8$-lattice. From the group law on the regular
fibers the distinguished singular nodal rational fiber minus the node inherits a
group law isomorphic to the multiplicative group $\C^\times$. Looking at
intersection points of the Mordell-Weil lattice with the distinguished 
singular fiber gives the isomorphism (\ref{eq_4.9}).

Combining the above two theorems gives a geometric interpretation of diagram
(\ref{eq_3.33}) for $R$ of type $E_8$ and $k = \frac{1}{6}$.

\begin{theorem}\label{theo_4.4}
We have a commutative diagram
\begin{equation*}
\begin{array}{ccc}
{\mathcal M}' & \longrightarrow & \Gamma'_8\backslash \mB^8\\
\downarrow && \downarrow\\
{\mathcal M} & \longrightarrow & \Gamma_8 \backslash \mB^8
\end{array}
\end{equation*}
with the bottom arrow given by the top arrow of the diagram 
in Theorem (\ref{theo_4.2}) and the
top arrow being the bottom arrow of diagram (\ref{eq_3.33}) with $R$ of type
$E_8$ and $k = \frac{1}{6}$ in accordance with diagram (\ref{eq_4.9}).
\end{theorem}

\begin{remark}\label{rem_4.5}
For $R$ of type $E_8$ with Coxeter number $h = 30$, the exponent difference along
the blown up identity element equals $1-hk = 1-30/6 = -4$. Since $w_0 = -1$ the
exponent difference on the quotient by $W$ of this divisor becomes $-2$ which is
not the reciprocal of an integer. Apparently the Schwarz conditions only hold
after exploitation of the hidden symmetry on $W\backslash H^\circ$. Since
$\Gamma'_8$ need not be a normal subgroup of $\Gamma_8$ this hidden symmetry is
just an equivalence relation on $W\backslash H^\circ$ with classes of 12 element
sets, and need not be obtained from an action of a group of order 12 on
$W\backslash H^\circ$. A description of this hidden symmetry in $E_8$ root system
terms remains unclear.
\end{remark}


\begin{thebibliography}{99}

\bibitem{ACT02} 
D.~Allcock, J.~Carlson and D.~Toledo: 
{\it The complex hyperbolic geometry of the moduli space of cubic surfaces}, 
J.\ Alg.\ Geom.\ 11 (2002), 659--724.

\bibitem{BD96} 
A.~Beilinson and V.~Drinfeld:
{\it Quantization of Hitchin's fibration and Langlands program}, 
Math.\ Phys.\ Studies 19 (1996), 3--7.


\bibitem{BD97} 
A.~Beilinson and V.~Drinfeld: 
{\it Quantization of Hitchin's integrable system and Hecke eigensheaves}, 
preprint (1997).

\bibitem{BN03} 
D.~Ben-Zvi and T.~Nevins:
{\it From solitons to many body systems}, \texttt{math.AG/0310490}.

\bibitem{B68} 
N.~Bourbaki:
{\it Groupes et alg\`{e}bres de Lie, Chapitres 4, 5 et 6},
Masson, Paris 1968.

\bibitem{C91} 
I. Cherednik: 
{\it A unification of Knizhnik-Zamolodchikov equations and
Dunkl operators via affine Hecke algebras}, 
Invent.\ Math.\ 106 (1991), 411--432.

\bibitem{CHL03} 
W.~Couwenberg, G.~Heckman and E.~Looijenga: 
{\it Geometric structures on the complement of a projective arrangement}, 
\texttt{math.AG/0311404}.

\bibitem{DM86} 
P.~Deligne and D.~Mostow: 
{\it Monodromy of hypergeometric functions and
non-lattice integral monodromy}, 
Publ.\ Math.\ IHES 63 (1986), 58--89.

\bibitem{D86} 
V.G.~Drinfeld: 
{\it Degenerate affine Hecke algebras and Yangians}, 
Funct.\ Anal.\ Appl.\ 20 (1986), 58--60.

\bibitem{H91} 
G.J.\ Heckman:
{\it An elementary approach to the hypergeometric shift operators of Opdam}, 
Invent.\ Math.\ 103 (1991), 341--350.

\bibitem{H97} 
G.J.~Heckman: 
{\it Dunkl operators}, S\'eminaire Bourbaki 828, Ast\'erisque 245
(1997), 223--246.

\bibitem{HL02} 
G.~Heckman and E.~Looijenga: 
{\it The moduli space of rational elliptic surfaces}, 
Alg.\ Geom.\ 2000, Azumino, Adv. Stud. Pure Math. 36 (2002), 185--248.

\bibitem{HO87/88} 
G.J.~Heckman and E.M.~Opdam: 
{\it Root systems and hypergeometric
functions I} (by H-O), {\it II} (by H), {\it III} (by O), {\it IV} (by O), 
Comp.\ Math. 64 (1987), 329--352 and 353--373, and 67 (1988), 21--49 and 191--209.

\bibitem{HM99} J.C.~Hurtibise and E.~Markman: 
{\it Calogero-Moser systems and Hitchin systems}, 
\texttt{math. AG/9912161}.

\bibitem{K00} S.~Kondo: 
{\it A complex hyperbolic structure on the moduli space of
curves of genus three}, J.\ reine u.\ angew.\ Math. 525 (2000), 219--232.

\bibitem{L98} E. Looijenga:
{\it Affine Artin groups and the fundamental groups of some
moduli spaces}, \texttt{math. AG/9801117}.

\bibitem{M03} I.G.~Macdonald: 
{\it Affine Hecke algebras and orthogonal polynomials},
Cambridge Tracts in Math.\ 157, CUP, Cambridge 2003.

\bibitem{M92} A.~Matsuo: 
{\it Integrable connections related to zonal spherical
functions}, Invent.\ Math.\ 110 (1992), 95--121.

\bibitem{M75} J. Moser: 
{\it Three integrable Hamiltonian systems connected with
isospectral deformation}, Adv.\ in Math. 16 (1975), 197--220.

\bibitem{OP76} M.A.~Olshanetsky and A.M.~Perelomov: 
{\it Completely integrable Hamiltonian systems connected with semisimple 
Lie algebras}, Invent.\ Math.\ 37 (1976), 93--108.

\bibitem{OP81} M.A.~Olshanetsky and A.M.~Perelomov: 
{\it Classical integrable finite dimensional systems related to Lie algebras}, 
Physics Reports 71, No.\ 5 (1981),
313--400.

\bibitem{OP83} M.A.~Olshanetsky and A.M.~Perelomov: 
{\it Quantum integrable systems
related to Lie algebras}, Physics Reports 91, No.\ 6 (1983), 314--403.

\bibitem{O95} E.M.~Opdam: 
{\it Harmonic analysis for certain representations of graded
Hecke algebras}, Acta Math.\ 175 (1995), 75--121.

\bibitem{O00} E.M.~Opdam: 
{\it Lecture notes on Dunkl operators for real and complex
reflection groups}, Math.\ Soc.\ Japan Memoirs, Vol.\ 8, 2000.

\bibitem{V99} R.~Vakil: 
{\it Twelve points on the projective line, branched covers, and
rational elliptic fibrations}, 
Math.\ Ann.\ 320 (2001), no. 1, 33--54.

\end{thebibliography}
\end{document}